\documentclass[letter, amsfonts, amssymb, reprint, showkeys, nofootinbib, twoside]{revtex4-1}
\usepackage[english]{babel}
\usepackage[utf8]{inputenc}
\usepackage[colorinlistoftodos, color=green!40, prependcaption]{todonotes}
\usepackage{graphicx} % Include figure files
\usepackage{dcolumn} % Align table columns on decimal point
\usepackage{bm} % bold math
\usepackage[intlimits]{amsmath}
\usepackage{esint}
\usepackage{caption}
\usepackage{verbatim}
\begin{document}

     %-------------------------------------------------------------------------
     % The introductory (header) part of the paper
     %-------------------------------------------------------------------------

     % The title of the paper. Use \shorttitle to indicate an abbreviated title
     % for use in running heads (you will need to uncomment it).

\title{An implementation of an efficient direct Fourier transform of polygonal areas and volumes}

\author{Brian B. Maranville}
    \email[Correspondence email address: ]{brian.maranville@nist.gov}% Your name
    \affiliation{NIST Center for Neutron Research \\
 100 Bureau Drive, Gaithersburg MD 20899 USA}

\date{\today} % Leave empty to omit a date

% \begin{synopsis}
% Calculations of the Fourier transform of a constant quantity over an area or volume defined by 
% polygons using the Divergence theorem to reduce the integral to a sum over vertices.
% \end{synopsis}

\begin{abstract}
Calculations of the Fourier transform of a constant quantity over an area or volume defined by 
polygons (connected vertices) are often useful in modeling wave scattering, or in fourier-space
filtering of real-space vector-based volumes and area projections.  If the system is 
discretized onto a regular array, Fast Fourier techniques can speed up the resulting calculations
but if high spatial resolution is required the initial step of discretization can limit performance;
at other times the discretized methods result in unacceptable artifacts in the resulting transform.
An alternative approach is to calculate the  
full Fourier integral transform of a polygonal area as a sum over the vertices, which has previously 
been derived in the literature using the divergence theorem to  reduce the problem from a 3-dimensional 
to line integrals over the perimeter of the polygon surface elements,
and converted to a sum over the straight segments of that contour.
We demonstrate a software implementation of this algorithm and show that it can provide accurate 
approximations of the Fourier transform of real shapes with faster convergence than a 
block-based (voxel) discretization.
\end{abstract}


\maketitle

\section {Introduction} \label{intro}
In calculations of scattering of neutrons, electrons, x-rays, etc. 
in the weak limit (the first Born approximation) it is useful
to be able to quickly calculate the Fourier transform (FT) of a constant quantity
such as the scattering power of a material over a well-defined volume.

The two most frequently used approaches to this problem are to use a combination of known
analytical transforms of shapes to create models \cite{Kline:do5025, SANS_model_functions}
or to discretize the volume on a regular, rectangular grid.
Both these strategies offer great speed in the calculation, as the analytic solution is typically a
small number of functions added together, while the gridded volume approach can take advantage of the 
inherent efficiency of Fast Fourier (FFT) methods.

There are however cases where it is not practical to construct a model from simple shapes
(the union of spheres example in this paper); 
at other times the discretization itself onto a rectangular 
grid is time-consuming or error-prone.  If a shape of interest has surfaces that are either
very curved or slightly misaligned with the coordinate axes, the gridded representation 
will contain large flat surface segments along the one of the underlying coordinates, 
resulting in artifacts in the FT that can dominate the calculation.  

It was first shown by Laue \cite{laue_polygon} and later expanded by others \cite{james1948optical} 
that using the divergence theorem, it is possible to directly calculate the 
Fourier transform of a volume element defined by a polygonal surface mesh.  Now this type of 
surface parametrization is regularly used in 3-dimensional graphics (openGL) as well as
in finite-element modeling programs such as the micromagnetic modeling program 
NMAG \cite{nmag_ieee}, the surface-meshing program netgen \cite{netgen_article, netgen_program}, 
and other pde solvers.

\section{Divergence theorem method in 2 dimensions}
As was shown earlier in \cite{laue_polygon, james1948optical}, 
one can approach the Fourier integral of polygonal volumes by 
beginning with the divergence theorem in 2 dimensions:
\begin{equation}
    \label{eq:divergence}
    \iint_A\left(\nabla\cdot\mathbf{F}\right)dA=\oint_C \mathbf{F} \cdot \mathbf{\hat n} \, ds
\end{equation}
where $\mathbf F$ is an arbitrary function of $(x,y)$, $C$ is the perimeter of the area 
(traversed counter-clockwise,) and $\mathbf{\hat n}$ is the unit vector normal to that perimeter 
and pointing outward from the defined area. 
This can be used to evaluate the FT$_A$ integral by constructing $\mathbf F$
as follows\cite{kazhdan2005reconstruction}:
\begin{equation}
    \mathbf{F} = 
    \begin{pmatrix}
        \frac{-i Q_x}{Q_x^2 + Q_y^2} e^{i(Q_x x + Q_y y)} \\
        \frac{-i Q_y}{Q_x^2 + Q_y^2} e^{i(Q_x x + Q_y y)}
    \end{pmatrix}
\end{equation}

One can verify that $\nabla \cdot \mathbf{F} = e^{i(Q_x x + Q_y y)}$, which is the original integrand.
Then for a perimeter line segment $S$ (again, traversing $C$ counter-clockwise) with 
start point 
$\mathbf p_0 \equiv \begin{pmatrix} x_0 \\ y_0 \end{pmatrix}$ and end point
$\mathbf p_1 \equiv \begin{pmatrix} x_1 \\ y_1 \end{pmatrix}$, the normal is
\begin{equation}
    \label{eq:line_normal}
    \mathbf{\hat n}_S = 
    \begin{pmatrix}
        (y_1-y_0)/d \\
        (x_0-x_1)/d
    \end{pmatrix}
\end{equation}
where $d=\sqrt{(x_1-x_0)^2 + (y_1-y_0)^2}$ is the length of the segment.
The RHS of Eq. \ref{eq:divergence} for that segment becomes
\begin{widetext}
\begin{equation}
\begin{array}{c}
    \label{eq:line_integral_eval}
    \int_{\mathbf{r}_0}^{\mathbf{r}_1} \mathbf{F} \cdot \mathbf{\hat n} \, ds = 
      \frac{-i\left( Q_x\frac{y_1-y_0}{d} + Q_y\frac{x_0 - x_1}{d} \right)}{Q_x^2 + Q_y^2} 
      \int_{\xi=0}^d e^{i(Q_x [x_0 + (x_1-x_0)\xi/d] + Q_y[y_0 + (y_1-y_0)\xi/d]} d\xi \\ [1.5em]
    = \frac{1}{Q_x^2 + Q_y^2} 
      \frac{-i\left( Q_x\frac{y_1-y_0}{d} + Q_y\frac{x_0 - x_1}{d} \right)}{i\left( Q_x\frac{x_1-x_0}{d} + Q_y\frac{y_1-y_0}{d} \right)}
      \left[ e^{i(Q_x x_1 + Q_y y_1)} - e^{i(Q_x x_0 + Q_y y_0)} \right] \\ [1.5em]
    = \frac{1}{Q_x^2 + Q_y^2} 
      \frac{ -Q_x \Delta_y + Q_y \Delta_x}{ Q_x \Delta_x + Q_y \Delta_y}
      \left[ e^{i(Q_x x_1 + Q_y y_1)} - e^{i(Q_x x_0 + Q_y y_0)} \right] \\
\end{array}
\end{equation}
\end{widetext}
where $\Delta_x =x_1 - x_0, \, \Delta_y = y_1-y_0$.  
We can rewrite Eq. \ref{eq:line_integral_eval} more generically
for any segment
\begin{equation}
    \label{eq:line_integral_generic}
    \displaystyle
    \int\limits_{\mathbf{r}_l}^{\mathbf{r}_{l+1}} \mathbf{F} \cdot \mathbf{\hat n} \, ds = 
      \frac{-(\mathbf{Q}\cdot \hat n) |\mathbf{d}_l|}{Q^2}
      \left(\frac{e^{i\mathbf{Q}\cdot\mathbf{r}_{l+1}}-e^{i\mathbf{Q}\cdot\mathbf{r}_{l}}}
      {\mathbf{Q}\cdot\mathbf{d}_l} \right)
\end{equation}
Where $\mathbf{d}_l \equiv \mathbf{r}_{l+1} - \mathbf{r}_l$. 
Factoring out $e^{i\mathbf{Q}\cdot(\frac{\mathbf{r}_{l+1} - \mathbf{r}_{l}}{2})}$ and noting
that $\hat n$ is defined to be perpendicular to $\mathbf{d}_l$, one gets
\begin{equation}
  \begin{array}{l l}  
    \int\limits_{\mathbf{r}_l}^{\mathbf{r}_{l+1}} 
    \mathbf{F} \cdot \mathbf{\hat n} \, ds &= 
    \frac{- Q d_l \sin\theta}{Q^2}
    e^{i \mathbf{Q} \cdot (\frac{\mathbf{r}_{l+1} + \mathbf{r}_{l}}{2})} 
    \frac{i \sin(\frac{Q d_l \cos\theta}{2})}{\frac{Q d_l \cos\theta}{2}} \\
    {} &= \frac{-2i\tan\theta }{Q^2} 
    e^{i \mathbf{Q} \cdot (\frac{\mathbf{r}_{l+1} + \mathbf{r}_{l}}{2})}
    \sin(\frac{Q d_l \cos\theta}{2})
  \end{array}    
\end{equation}
where $\theta$ is the angle between $\mathbf{Q}$ and $\mathbf{d}_l$.

It is interesting to note that for a closed polygon the exponential $e^{i(Q_x x_n + Q_y y_n)}$ will show up exactly
twice in the sum, once as the $(x_0, y_0)$ term and once as the $(x_1, y_1)$ term.  Thus the sum can be recast as
a sum over all the vertices instead of all the segments, with a weighting that depends on the difference of the 
tangent of the angle between $\mathbf{Q}$ and the segment for the two segments meeting at that point.

The form of Eq. \ref{eq:line_integral_generic} is preferred though, because in the limit of 
$\mathbf{Q}\cdot(\mathbf{r}_{l+1}-\mathbf{r}_l) \rightarrow 0$ the quantity in parentheses at the right will not diverge.

\section{Divergence theorem method in 3 dimensions}

Using the same logic to extend to three dimensions, the integral becomes
\begin{equation}
    \mathrm{FT}_V \equiv \iiint_V  e^{i\mathbf{Q}\cdot\mathbf{r}} \,dV
\end{equation}
over a volume $V$ with a polygonal surface (made up of connected planar polygons) such as a tetrahedron or cube.
The three-dimensional divergence theorem can be stated as 
\begin{equation}
    \label{eq:divergence3d}
    \iiint_V\left(\nabla\cdot\mathbf{F}\right)dV=\oiint_S \mathbf{F} \cdot \mathbf{\hat n_S} \, dS
\end{equation}
Because there are a finite number of planar polygons that make up the surface, the integral
can be decomposed to a sum of integrals over those regions.
\begin{eqnarray}
  \displaystyle
  %\begin{array}{c}
    \oiint_S \mathbf{F} \cdot \mathbf{\hat n_S} \, dS = \sum_m I_m \\ [1.5em]
    I_m = \iint_{A_m} \mathbf{F}\cdot\mathbf{\hat n}_m dA
  %\end{array}
\end{eqnarray}
Where $m$ is the number of polygons on the surface, and $\mathbf{\hat n}_m$ is the surface normal for area $A_m$.

Again we define
\begin{equation}
    \mathbf{F} = \frac{-i \mathbf{Q}}{Q^2} e^{i\mathbf{Q}\cdot\mathbf{r}}
\end{equation}
and again, $\nabla \cdot \mathbf{F} = e^{i\mathbf{Q} \cdot \mathbf{r}}$, the original integrand.
For a given surface polygon, 
\begin{equation}
    \iint_A \mathbf{F}\cdot\mathbf{\hat n} \, dA = \iint_{A_\perp \hat n} \frac{-iQ_\parallel}{Q^2} e^{i\mathbf{Q}\cdot\mathbf{r}} \,dA
\end{equation}
Where $Q_\parallel$ is the component of $\mathbf{Q}$ along $\mathbf{\hat n}$.  
Since the area being integrated over is by definition perpendicular to the surface normal, we can pull out that
component of the integrand as a constant
\begin{equation}
    \label{eq:single_area_integral}
    \iint_A \mathbf{F}\cdot\mathbf{\hat n} \, dA = \frac{-iQ_\parallel}{Q^2} e^{iQ_\parallel r_\parallel} \iint_{A_\perp \hat n}  e^{i\mathbf{Q}_\perp\cdot\mathbf{r}_\perp} dA
\end{equation}
where $\mathbf{Q}_\perp , \mathbf{r}_\perp$ are the components perpendicular to $\mathbf{\hat n}$ (in the integration plane.)

Now the volume integral is reduced to a sum of area integrals of the same type as was treated in the previous section.
Inserting Eq. \ref{eq:line_integral_generic} into Eq. \ref{eq:single_area_integral},
we get an expression like the one found in Eq. 10.80 of \cite{james1948optical}:
\begin{widetext}
\begin{equation}
  \label{eq:sum_elements_3d}
    %\iiint_V  e^{i\mathbf{Q}\cdot\mathbf{r}} \,dV = 
    \mathrm{FT}_V = 
    \sum_m I_m = \sum_m
      \frac{-iQ_{\parallel m}}{Q^2} e^{iQ_{\parallel m} r_{\parallel m}} \sum_{l=0}^{\nu_m + 1} 
      \frac{-\mathbf{Q}_{\perp m}\cdot(\mathbf{\hat n}_m \times(\mathbf{r}_{l+1}-\mathbf{r}_l))}{Q_{\perp m}^2}
      \left(\frac{e^{i\mathbf{Q}_{\perp m} \cdot\mathbf{r}_{l+1}}-e^{i\mathbf{Q}_{\perp m} \cdot\mathbf{r}_{l}}}{\mathbf{Q}_{\perp m} \cdot(\mathbf{r}_{l+1}-\mathbf{r}_l)} \right)
\end{equation}
\end{widetext}
where the last sum is a counterclockwise sum over vertices $\mathbf{r}_l$ in the contour $C_m$ defining area element $A_m$, and
the notation for counting vertices along $C_m$ is understood to wrap around, so that if there are $\nu_m$ vertices
defining $C_m$ then $\mathbf{r}_{\nu_m +1} \equiv \mathbf{r}_0$.
Also, ${}_{\parallel m}, {}_{\perp m}$ refer to components parallel and perpendicular to the the surface normal $\mathbf{\hat{n}}_m$
of area element $A_m$.

\section{Application to real geometries}
\subsection{Transform of a rectangle}
\begin{figure}
	\includegraphics[width=\linewidth]{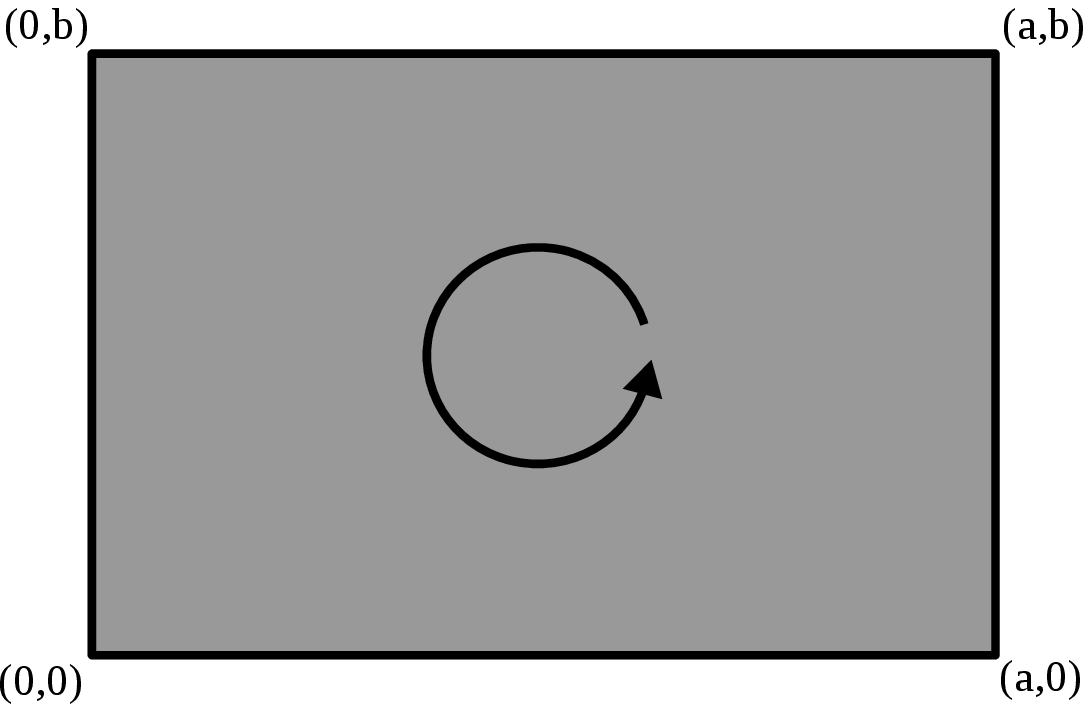}
	\caption{Rectangle described in the example, with lower-left corner at origin with width $a$ and height $b$.  The arrow in the center indicates the direction of traversal (counterclockwise) when doing the sum over the segments and vertices.}
	\label{fig:rectangle_diag}
\end{figure}

As an example of a 2-d system, let's take a rectangle of length $a$ along the $x$ axis and $b$ along the $y$ axis, as 
seen in Fig. \ref{fig:rectangle_diag}.
The definite integral in this case is separable in $x$ and $y$, and we get 
\begin{eqnarray}
    \label{eq:rectangle_integral}
    \iint\limits_{\mathrm{Rectangle}} e^{i\mathbf{Q}\cdot\mathbf{r}} \, dA &=& 
    \int_{0}^{a} e^{iQ_x x} dx \int_{0}^{b} e^{iQ_y y} dy \\ \nonumber
    &=& \frac{-1}{Q_x Q_y}(e^{iQ_x a} - 1)(e^{iQ_y b} - 1)
\end{eqnarray}

If we travel counterclockwise around the rectangle we get these vertices 
\begin{equation}
  \begin{array}{c}
    \mathbf{r}_0 = (0,0) \\
    \mathbf{r}_1 = (a,0) \\
    \mathbf{r}_2 = (a,b) \\
    \mathbf{r}_3 = (0,b) 
  \end{array}
\end{equation}
and putting them into Eq. \ref{eq:line_integral_eval}, 
\begin{equation}
  \label{eq:rectangle_sum}
  \begin{array}{lll}
    \oint_C \mathbf{F} \cdot \mathbf{\hat n} \, ds =&\,\,&
      \frac{1}{Q^2} \frac{ Q_y a}{ Q_x a} \left[ e^{iQ_x a} - e^{i0} \right] \\
   &+&\frac{1}{Q^2} \frac{ -Q_x b}{ Q_y b} \left[ e^{i(Q_x a + Q_y b)} - e^{iQ_x a} \right] \\
   &+&\frac{1}{Q^2} \frac{ -Q_y a}{ -Q_x a} \left[ e^{iQ_y b} - e^{i(Q_x a + Q_y b)} \right] \\
   &+&\frac{1}{Q^2} \frac{ Q_x b}{ -Q_y b} \left[ e^{i0} - e^{iQ_y b} \right] \\
  \end{array}
\end{equation}
and a bit of algebra yields 
\begin{equation}
    FT_{\mathrm{Rect}} = \frac{1}{Q^2} \frac{Q_x^2 + Q_y^2}{Q_x Q_y}\left[e^{iQ_y b} + e^{iQ_x a} - e^{i(Q_x a + Q_y b)} - 1\right] 
\end{equation}
which is equivalent to Eq. \ref{eq:rectangle_integral}

\subsection{Rectangular prism}
Now extending to the third dimension, we take a rectangular prism of length $a$ along the $x$ axis, 
$b$ along the $y$ axis and $c$ along the $z$ axis.  The definite integral is then
\begin{equation}
	\label{eq:prism_direct}
	\iiint\limits_\mathrm{Prism} e^{i\mathbf Q \cdot \mathbf r} \, dV =
		\frac{i(e^{i Q_x a} - 1)(e^{i Q_y b} - 1)(e^{i Q_z c} - 1)}{Q_x Q_y Q_z}
\end{equation}

Breaking it up into the 6 surfaces that compose the boundary of this region, we get from 
Eq. \ref{eq:single_area_integral} the following for the top surface $(z=c)$:
\begin{equation}
	I_{z=c} = \frac{i Q_z e^{i Q_z c}}{Q^2} 
	  \frac{(e^{i Q_x a} - 1)(e^{i Q_y b} - 1)}{Q_x Q_y}
\end{equation}
while for the bottom surface $(z=0)$ we get
\begin{equation}
	I_{z=0} = \frac{-i Q_z}{Q^2} 
	  \frac{(e^{i Q_x a} - 1)(e^{i Q_y b} - 1)}{Q_x Q_y}
\end{equation}
Adding those two pieces (and multiplying by $Q_z$ on the top and bottom) we get
\begin{equation}
	I_{z=c} + I_{z=0} = \frac{i Q_z^2}{Q^2} 
	  \frac{(e^{i Q_x a} - 1)(e^{i Q_y b} - 1)(e^{i Q_z c} - 1)}{Q_x Q_y Q_z}
\end{equation}
and clearly the sum of all six interfaces will give the same answer as in Eq. \ref{eq:prism_direct} above.

\subsection{Sphere}
The Fourier transform of a sphere is known analytically as well, so we can compare the 
discretization methods for a model that does not naturally align with a rectilinear grid.
The analytic solution for the unit sphere in 3-d is 
\begin{equation}
	\mathrm{FT}_{V}(Q) =\frac{4\pi}{3} \Gamma(5/2) \frac{\mathrm{J}_{3/2}(Q)}{(Q/2)^{3/2}} ~ (\mathrm{sphere})
\end{equation}
where $\Gamma$ is the gamma function, and J$_{3/2}$ is a half-integer Bessel function of the first kind.

\begin{figure}
	\includegraphics[width=\linewidth]{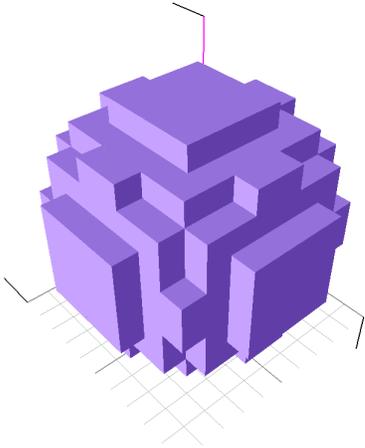}
	\caption{Voxelization of the unit sphere, with voxel dimension = 0.2.
	This results in 522 non-zero voxels}
	\label{fig:sphere_voxels}
\end{figure}

Also, the sphere was voxelized using the open-source program \emph{binvox}\cite{10.1109/TVCG.2003.10002}\cite{binvox}
and the results can be seen in Fig. \ref{fig:sphere_voxels}. 
% To make the VRML .wrl file, choose VRML as mesh output type in netgen, 
%   then export mesh as "sphere.wrl"
% to make the voxelization:
% ./binvox -d 10 ~/design/papers/bb_manuscripts/fourier_polygon/sphere.wrl
% to view the voxelization:
% ./viewvox ~/design/papers/bb_manuscripts/fourier_polygon/sphere_5.binvox 

\begin{figure}
	\includegraphics[width=\linewidth]{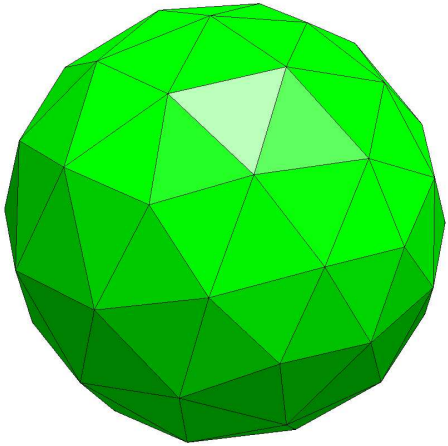}
	\caption{Surface meshing of the unit sphere.  The mesh consists of 120
	triangles connecting 83 nodes.}
	\label{fig:sphere_polygon}
\end{figure}

Finally, a surface mesh as rendered by the open-source meshing program
\emph{netgen}\cite{netgen_article}\cite{netgen_program} is seen in Fig. \ref{fig:sphere_polygon}.
The meshing input radius was adjusted so that the volume of the resulting polygon matches that of the
unit sphere.

\begin{figure}
	\includegraphics[width=\linewidth]{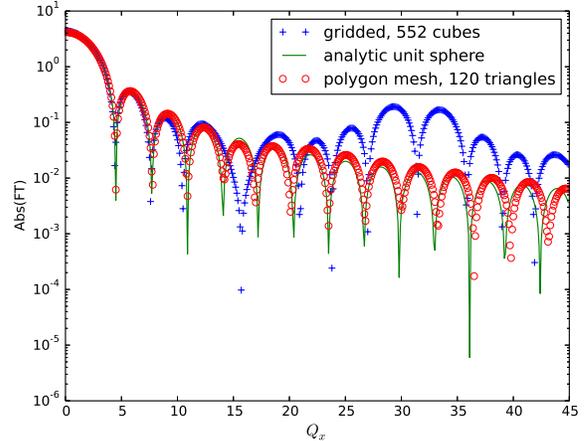}
	\caption{FT$_V$ of the unit sphere (along $Q_x$) as calculated 
	analytically (solid green line), with the grid voxelization (blue plus symbols)
	and with the polygon approximation to the surface (red open circles)}
	\label{fig:sphere_FT_compare}
\end{figure}

The calculation of FT from each of these methods is presented in Fig \ref{fig:sphere_FT_compare}, 
calculated along the $Q_x$-axis (the artifacts from gridding are likely to be most apparent along
the coordinate axes.)
We do see artifacts in the gridded-FT data, as expected at $Q_x \approx 31.4 \approx 2\pi/d$ where $d=0.2$ is the 
discretization size of the grid.  The polygon-FT appears to be a good approximation to the analytic
solution even out to high $Q_x$, and a much better approximation than the gridded FT.

\subsection{Union of spheres}
\label{sec:spheres_example}
\begin{figure}
	\includegraphics[width=\linewidth]{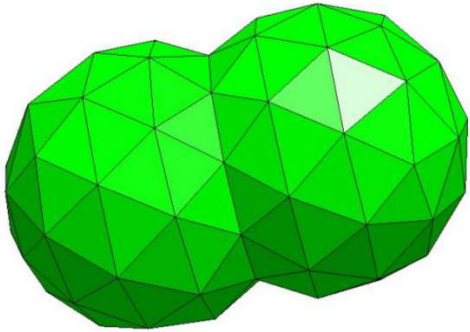}
	\caption{Surface meshing of the volume created by the union of two spheres.
	The mesh consists of 160 
	triangles connecting 83 nodes.  No refinement of the mesh was performed beyond the 
	initial meshing.}
	\label{fig:2spheres_mesh}
\end{figure}

Consider a volume defined by the union of two spheres with radius $r = 1.0$ 
and centers at $x=-0.6$ and $x=0.6$.  
This is an example of a volume for which it is difficult to calculate the FT 
analytically, as it would require subtracting the FT of the lens-shaped overlap from the
sum of the two sphere FT functions. 
A surface mesh of this geometry 
can be seen in Fig. \ref{fig:2spheres_mesh} (again rendered by 
%the open-source meshing program
\emph{netgen}.)

\begin{figure}
	\caption{$Q_x,Q_y$ projection of the absolute value of the Fourier transform of the union of two spheres.  
	(Calculated at $Q_z=0$)}
	\includegraphics[width=\linewidth]{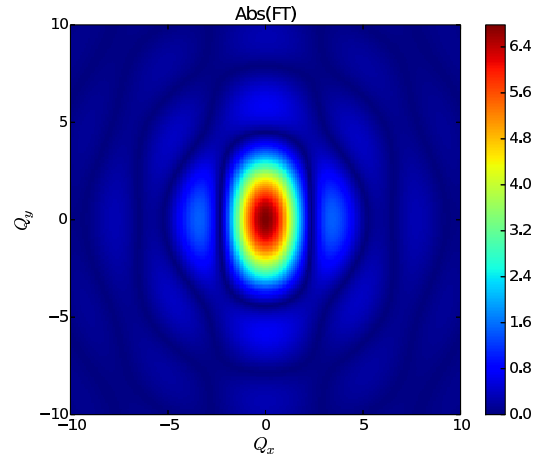}
	
	\label{fig:2spheres_absFT}
\end{figure}

The absolute value of FT for this polygonal shell in the $Q_x,Q_y$ plane is plotted
as a heatmap in Fig. \ref{fig:2spheres_absFT}, with the color scale on the right.
%, and the imaginary component is plotted in Fig. \ref{fig:2spheres_imagFT}.  The geometry that is 
%described is centro-symmetric, and so the FT is expected to be purely real-valued; in this case
%the magnitude of the imaginary component corresponds to deviations of from the analytical solution to
%the FT of the union of two ideal spheres.  The fact that the
%imaginary component is so small compared to the absolute value indicates that the relatively coarse
%meshing used already provides a fairly good approximation to the geometry as far as the FT is concerned.

%\begin{center}
%	\includegraphics[width=\linewidth]{imagFT}
%	\captionof{figure}{$xy$-projection of the imaginary component of the Fourier transform of the union of two spheres.  
%	(Calculated at $Q_z=0$)}
%	\label{fig:2spheres_imagFT}
%\end{center}

\section{Conclusions}
We have shown that by using a direct calculation of the Fourier transform of a constant function
over a volume defined by a connected-polygon surface, a reasonable approximation to the FT of a 
sphere can be achieved from a discretization to just 120 surface triangles.  This is compared to a 
approximation of the same shape using a much larger number of regularly-gridded volume elements
(voxels).  The same technique is shown to be straightforward to apply to complicated shapes
that are hard to calculate the analytic FT for directly, such as the overlapping spheres.
Artifacts from the conversion to a rectilinear basis are avoided, and this should provide a more
straightforward method of calculating the FT from volumes where the surface mesh is already known, 
such as can be extracted from tomography, microscopy or other real-space probes, without requiring
the additional (computationally intensive) process of voxelization.

% \ack{Acknowledgements}
\begin{acknowledgments}
The author would like to thank Prof. M. Hore of CWRU and Drs. B. Kirby and K. Krycka of NIST
for helpful discussions.
\end{acknowledgments}

% \bibliographystyle{revtex4-1} 
% \bibliography{master}
%merlin.mbs apsrev4-1.bst 2010-07-25 4.21a (PWD, AO, DPC) hacked
%Control: key (0)
%Control: author (72) initials jnrlst
%Control: editor formatted (1) identically to author
%Control: production of article title (-1) disabled
%Control: page (0) single
%Control: year (1) truncated
%Control: production of eprint (0) enabled
%

\newpage
\clearpage

\appendix
\begin{widetext}

\section{Support files and code}
\subsection{Geometry definition for union of two spheres}
Using the geometry file format for \emph{netgen}:
\begin{verbatim}
algebraic3d

solid main = sphere (-0.6, 0, 0; 1);
solid second = sphere (0.6, 0, 0; 1);

solid combined = main or second;
tlo combined;
\end{verbatim}

\subsection{Python code for calculations}
This includes a function for reading in a surface mesh file that can be 
written from \emph{netgen}, as well as a demo of the main function 
for calculating the Fourier transform of the example in section \ref{sec:spheres_example}.
%\lstinputlisting[language=Python]{fourier_surfacemesh.py}
%\verbatiminput{fourier_surfacemesh.py}
\begin{verbatim}
import numpy as np
eps = 1e-30

def read_surface_meshfile(fileobj):
    f = fileobj
    if f.readline().strip() != "surfacemesh":
        print("not a surface mesh file")
    numpoints = int(f.readline().strip())
    points = np.empty((numpoints,3), dtype="float")
    for i in range(numpoints):
         points[i] = np.array(f.readline().strip().split(), dtype='float')
    numelements = int(f.readline().strip())
    elements = np.empty((numelements, 3), dtype='int')
    for j in range(numelements):
        elements[j] = np.array(f.readline().strip().split(), dtype='int')
    # indexing of points in the geometry file begins with 1, 
    # but it begins with 0 in numpy ndarray, so need to subtract 1 
    # from every point in element index array:
    return {'points': points, 'elements': elements - 1}
            
def resolve_coords(points, elements, wrap=True):
    # convert elements with point ids to lists of coords
    # if wrap is True, add point[0] to the end of each list
    el = elements.copy()
    if wrap:
        el = np.concatenate((el, el[:,:1]), axis=1)
    return (points[el]).copy()
    
def get_normal_vec(p):
    v1 = p[:,1] - p[:,0]
    v2 = p[:,2] - p[:,0]
    normal = np.cross(v1, v2)
    normal = normal * 1.0 / np.sqrt(np.sum(normal*normal, axis=1))[:,None] 
    return normal
    
def fourier_vec(qx, qy, qz, p):
    """ calculate the fourier transform of the volume bounded by the elements
    in the list p"""
    ########################################################################## 
    # Inputs:                                                                #
    # qx, qy and qz should have 3 dimensions each but they can be sparse,    #
    # i.e. qx.shape can be (4,1,1) when qy.shape is (1, 12, 1) etc.          #
    # in which case they will be broadcast.                                  #
    #                                                                        #
    # p is an array of elements, which are themselves an array of points     #
    # describing a counterclockwise trip around the border of the element    #
    # (as seen from the outside of the volume) where a point is the array    #
    # [x,y,z]                                                                #
    #                                                                        #
    ##########################################################################
    # Outputs:                                                               #
    # result is an array covering all qx, qy and qz of the fourier transform #
    # of the volume enclosed by the elements in p                            #
    ##########################################################################
    normal = get_normal_vec(p)
    # now do dot product
    dotx = qx[:,:,:,None] * normal[:,0]
    doty = qy[:,:,:,None] * normal[:,1]
    dotz = qz[:,:,:,None] * normal[:,2]
    Qn_length = dotx + doty + dotz
    Qnx = Qn_length * normal[None, None, None, :,0]
    Qny = Qn_length * normal[None, None, None, :,1]
    Qnz = Qn_length * normal[None, None, None, :,2]
    Qpx = qx[:,:,:,None] - Qnx
    Qpy = qy[:,:,:,None] - Qny
    Qpz = qz[:,:,:,None] - Qnz
    Qp = np.concatenate((Qpx[...,None], Qpy[...,None], Qpz[...,None]), axis=-1)
    
    Qsq = qx**2 + qy**2 + qz**2
    Qpsq = np.sum(Qp * Qp, axis=-1)
    
    ## Note: p[:,0] is the first point in the element (for all elements)    
    rn_length = np.sum(normal * p[:,0], axis=1)
    subsum = np.zeros_like(Qn_length, dtype="complex")
    
    ## Here is equation 17 from the attached publication for a single element
    # (including the for loop which sums over the vertices)
    result = (1j * Qn_length / Qsq[:,:,:,None]) \
        * np.exp(1j * Qn_length * rn_length[None, None, None, :])
    
    for i in range(p.shape[1]-1):
        v1  = p[:,i+1] - p[:,i]
        sub1 = np.sum(Qp * np.cross(normal, v1)[None, None, None, :], axis=-1)
        sub1 = sub1 / (Qpsq + eps)
        sub2 = (np.exp(1j*np.sum(Qp*p[None,None,None,:,i+1],axis=-1)+eps/2.0) \
         -np.exp(1j*np.sum(Qp*p[None, None, None,:,i], axis=-1)-eps/2.0))
        sub3 = 1.0 / (np.sum(Qp*(p[None, None, None, :, i+1] \
         - p[None, None, None, :, i]), axis=-1) + 1j * eps)
        subsum += (sub1 * sub2 * sub3)
        
    return result * subsum
    
def demo():
    from pylab import figure, xlabel, ylabel, title, imshow, colorbar, show
    import StringIO
    
    qx = np.linspace(-10, 10, 40) + eps
    qy = np.linspace(-10, 10, 40) + eps
    qz = np.linspace(-0, 0, 1) + eps
    qx,qy,qz = np.meshgrid(qx,qy,qz, sparse=True)
    extent = (qx.min(), qx.max(), qy.min(), qy.max())
    
    surf = read_surface_meshfile(StringIO.StringIO(twospheres_surfacemesh))
    pp = resolve_coords(surf['points'], surf['elements'])
    ft = fourier_vec(qx, qy, qz, pp)
    fts = np.sum(ft[:,:,0], axis=-1)
    fig1 = figure()
    
    imshow(abs(fts), extent=extent, aspect=1)
    title('Abs(FT)')
    xlabel('$Q_x$', size='large')
    ylabel('$Q_y$', size='large')
    colorbar()
    
    fig2 = figure()
    imshow(fts.imag, extent=extent, aspect=1)
    title('Imaginary FT')
    xlabel('$Q_x$', size='large')
    ylabel('$Q_y$', size='large')
    colorbar()    
    show()
    
   
twospheres_surfacemesh = """\
surfacemesh
83
         0       -0.8          0 
         0  -0.646633  -0.471026 
         0  -0.244492  -0.761724 
         0    0.24951  -0.760095 
         0   0.650833  -0.465206 
         0   0.799973 0.00660409 
         0   0.644364   0.474125 
         0   0.242392   0.762395 
         0  -0.250252   0.759851 
         0  -0.648316   0.468707 
 -0.356907  -0.922165  -0.300863 
 -0.347131  -0.567188  -0.783808 
 -0.357573  0.0050912  -0.970156 
 -0.369647   0.575911  -0.784388 
  -0.37288   0.931353  -0.284602 
  -0.37676   0.917429   0.329375 
 -0.369284   0.561424   0.794716 
 -0.387142 -0.0140678   0.976982 
 -0.412157  -0.582938   0.790505 
 -0.385578  -0.929431   0.300301 
 -0.760251  -0.306731  -0.938209 
 -0.760576  -0.795816   -0.58386 
 -0.813869   0.298987  -0.929982 
 -0.813426    0.79047  -0.574114 
 -0.821753   0.973827  0.0498692 
 -0.835588   0.746598   0.622165 
 -0.761439   0.283083   0.945411 
 -0.903893  -0.261241   0.916189 
 -0.883309  -0.779017   0.559346 
 -0.803356  -0.978812 -0.0239454 
  -1.19796 -0.0347093  -0.800775 
  -1.18082  -0.758163  -0.296374 
  -1.23593   0.456043  -0.622585 
  -1.20125   0.760906  -0.243974 
    -1.276   0.676061   0.293208 
  -1.22172   0.230192   0.748648 
  -1.10981  -0.519442  -0.685767 
   -1.1967  -0.780973   0.184463 
  -1.47888  0.0876796  -0.468915 
  -1.38991  -0.362765   0.494419 
  -1.50034   0.409012  -0.148656 
  -1.52632   0.184822   0.328273 
  -1.43381  -0.350022  -0.426898 
  -1.45293  -0.521864 -0.0129389 
  -1.59324  -0.111427 -0.0324898 
  0.353411  -0.922424    -0.2972 
   0.35295  -0.576474  -0.778873 
  0.359646 -0.00259245  -0.970682 
  0.365788   0.577679  -0.781941 
  0.354885   0.929085  -0.276983 
  0.350692   0.915547   0.315627 
  0.356403   0.564491   0.788676 
  0.360185 -0.0185137   0.970642 
  0.354104  -0.571601   0.782821 
  0.355115   -0.92195   0.300067 
  0.778156  -0.806435  -0.563846 
  0.777153  -0.329327  -0.927448 
  0.807472   0.308424  -0.928348 
  0.815371   0.808817  -0.547202 
  0.775977   0.983647  0.0383495 
  0.772101   0.779013   0.602926 
  0.781753   0.290957   0.939314 
  0.789647  -0.311115   0.931258 
  0.779065  -0.781899   0.597134 
  0.783485  -0.982953 -0.0116784 
   1.13571  -0.519537  -0.665658 
   1.21105 -0.0477391  -0.790149 
   1.21305   0.486838   -0.62222 
   1.19736   0.784525  -0.166363 
   1.13036   0.790032   0.307512 
   1.19753   0.460595   0.656358 
   1.15168 -0.0067317   0.834029 
   1.23508  -0.439877    0.63497 
    1.2143  -0.753463   -0.23437 
   1.18111  -0.771965   0.257649 
   1.47664   -0.31298  -0.365446 
    1.4629   0.448967  -0.232032 
   1.44925   0.479567   0.220866 
   1.48355  0.0124405   0.468166 
   1.47929   0.115775  -0.461995 
    1.4961  -0.418884   0.146737 
    1.5993  0.0340127 -0.0153453 
 0.0209034 -0.00574919 -0.0184872 
160
       2       1      11
       3       2      12
       4       3      13
       5       4      14
       6       5      15
       7       6      16
       8       7      17
       9       8      18
      10       9      19
       1      10      20
       3      12      13
       4      13      14
       5      14      15
       6      15      16
       7      16      17
       8      17      18
       9      18      19
      10      19      20
       1      20      11
       2      11      12
      13      12      21
      12      11      22
      14      13      23
      15      14      24
      16      15      25
      17      16      26
      18      17      27
      19      18      28
      20      19      29
      11      20      30
      30      20      29
      29      19      28
      28      18      27
      27      17      26
      26      16      25
      25      15      24
      24      14      23
      23      13      21
      22      11      30
      12      22      21
      23      21      31
      22      30      32
      24      23      33
      23      31      33
      25      24      34
      24      33      34
      26      25      35
      25      34      35
      27      26      36
      26      35      36
      28      27      36
      29      28      40
      21      22      37
      31      21      37
      22      32      37
      30      29      38
      32      30      38
      29      40      38
      28      36      40
      33      31      39
      76      82      81
      34      33      41
      33      39      41
      35      34      41
      36      35      42
      40      36      42
      35      41      42
      79      81      82
      39      31      43
      31      37      43
      37      32      43
      32      38      44
      44      38      40
      43      32      44
      41      39      45
      42      41      45
      39      43      45
      43      44      45
      40      42      45
      44      40      45
       1       2      46
       2       3      47
       3       4      48
       4       5      49
       5       6      50
       6       7      51
       7       8      52
       8       9      53
       9      10      54
      10       1      55
       2      47      46
       3      48      47
       4      49      48
       5      50      49
       6      51      50
       7      52      51
       8      53      52
       9      54      53
      10      55      54
       1      46      55
      46      47      56
      47      48      57
      48      49      58
      49      50      59
      50      51      60
      51      52      61
      52      53      62
      53      54      63
      54      55      64
      55      46      65
      65      46      56
      64      55      65
      63      54      64
      62      53      63
      61      52      62
      60      51      61
      59      50      60
      58      49      59
      57      48      58
      47      57      56
      56      57      66
      57      58      67
      58      59      68
      59      60      69
      60      61      70
      61      62      71
      62      63      72
      63      64      73
      65      56      74
      56      66      74
      64      65      75
      73      64      75
      65      74      75
      72      63      73
      71      62      72
      70      61      71
      69      60      70
      68      59      69
      67      58      68
      66      57      67
      74      66      76
      66      67      76
      68      69      77
      70      71      78
      72      73      79
      76      67      80
      80      67      68
      75      81      73
      78      71      79
      68      77      80
      77      69      78
      69      70      78
      71      72      79
      75      74      81
      79      73      81
      74      76      81
      77      78      82
      80      77      82
      78      79      82
      76      80      82
"""

if __name__ == '__main__':
    demo()   
\end{verbatim}
\end{widetext}

\end{document}